\documentclass{fpsac} 
\usepackage{graphicx,psfrag}

\usepackage{algorithm,algorithmic,vmargin}
\usepackage{amssymb,amsmath,amscd,amsthm}

\def\FW{\mathcal{W}}   

\newtheorem{theorem}{Theorem}[section]

\theoremstyle{definition}

\theoremstyle{remark}

\numberwithin{equation}{section}

\newtheorem{prop}[theorem]{Proposition}

\title[{Generating Functions For Kernels of Digraphs}]{Generating
Functions For Kernels of Digraphs\\ (Enumeration \& Asymptotics for Nim Games)}

\author{Cyril Banderier}
\thanks{Corresponding author: Cyril.Banderier@lipn.univ-paris13.fr}
\email{Cyril.Banderier@lipn.univ-paris13.fr}
\urladdr{http://www-lipn.univ-paris13.fr/$\sim$banderier}

\author{Jean-Marie Le Bars} 
\email{Jean-Marie.Lebars@info.unicaen.fr}
\urladdr{http://users.info-unicaen.fr/$\sim$lebars}

\author{Vlady Ravelomanana}
\email{Vlady.Ravelomanana@lipn.univ-paris13.fr}
\urladdr{http://www-lipn.univ-paris13.fr/$\sim$ravelomanana}

\address{Cyril Banderier \& Vlady Ravelomanana\newline
CNRS (UMR 7030)/Laboratoire d'Informatique de Paris Nord \newline
Institut Galil\'ee, Universit\'e Paris 13 \newline
99, avenue Jean-Baptiste Cl\'ement \\ \newline
93430 Villetaneuse (France) }

\address{Jean-Marie Le Bars\newline
Laboratoire GREYC\newline
Universit\'e de Caen\newline
14032 Caen (France)}

\keywords{generating functions, analytic combinatorics, kernels of graphs, Nim games}

\begin{document}

\begin{abstract}
In this article, we study directed graphs (digraphs) with a coloring
constraint due to Von Neumann and related to Nim-type games.
This is equivalent to the notion of kernels of digraphs, which appears in numerous fields
of research such as game theory, complexity theory, 
artificial intelligence (default logic, argumentation in multi-agent systems), 
0-1 laws in monadic second order logic, combinatorics (perfect graphs)...
Kernels of digraphs lead to numerous difficult questions
(in the sense of NP-completeness, \#P-completeness).
However, we show here that it is possible to use a generating function
approach to get new informations: 
we use technique of symbolic and analytic combinatorics
(generating functions and their singularities) in order to get
exact and asymptotic results, e.g. for the existence of a kernel
in a circuit or in a unicircuit digraph. This is a first step toward
a generatingfunctionology treatment of kernels, while using, e.g., 
an approach ``\`a la Wright''.
Our method could be applied to more general ``local coloring constraints''
in decomposable combinatorial structures.\\
{\sc R\'esum\'e.} 
Nous \'etudions dans cet article les graphes dirig\'es (digraphes) avec une
contrainte de coloriage introduite par Von Neumann et reli\'ee aux jeux
de type Nim. Elle \'equivaut \`a la notion de noyaux de digraphes, qui
appara\^{\i}t dans de nombreux domaines, tels la th\'eorie des jeux, la th\'eorie de la
complexit\'e, l'intelligence artificielle (logique des d\'efauts,
argumentation dans les syst\`emes multi-agents), les lois 0-1 en logique
monadique du second ordre, la combinatoire (graphes parfaits)...
Les noyaux des digraphes posent de nombreuses questions difficiles 
(au sens de la NP-compl\'etude ou de la \#P-compl\'etude). Cependant,
nous montrons ici qu'il est possible de recourir aux s\'eries g\'en\'eratrices
afin d'obtenir de nouvelles informations~: nous utilisons les techniques
de la combinatoire symbolique et analytique (\'etude des singularit\'es
d'une s\'erie) afin d'obtenir des r\'esultats exacts ou asymptotiques, par
exemple pour l'existence d'un noyau dans un digraphe unicircuit. Il
s'agit l\`a de la premi\`ere \'etape vers une s\'erie g\'en\'eratrilogie des
noyaux. Notre m\'ethode peut \^etre appliqu\'ee plus g\'en\'eralement \`a des
``contraintes locales'' de coloriage dans des structures combinatoires
d\'ecomposables.
\end{abstract}

\maketitle

\section{Introduction}\label{Sec1}

Let $V$ and $E$ be the set of vertices and directed edges 
(also called {\em arcs}) 
of a directed graph $D$ without loops or multiarcs (we call such graphs {\em digraphs} hereafter).
A kernel of $D$ is a nonempty subset $K$ of $V$, 
such that for any $a, b\in K$, the edge $(a,b)$ does not belong to $E$, 
and for any vertex outside the kernel ($a \not \in K$), 
there is a vertex in the kernel ($b \in K$),  such that the edge $(a, b)$ belongs to $E$.
In other words, $K$ is a nonempty independent and dominating set of vertices in $D$~\cite{ber}.
Not every digraph has a kernel and if a digraph has a kernel, this kernel is not necessarily unique.
The notion of kernel allows elegant interpretations in various 
contexts, since it is related to other well-known concepts from graph theory  and complexity theory. 
In game theory the existence of a kernel corresponds to a
winning strategy in two players for famous Nim-type games (cf.~\cite{bs,fra,jap1,NeMo80}). 

Imagine that two players ${\mathcal A}$ and ${\mathcal B}$ play the following 
game on $D$ in which they move a token each in turn: ${\mathcal A}$ starts the game by choosing 
an initial vertex $v_0 \in V$ and then makes a move to a vertex $v_1$.
 A move consists in taking the token from the present position $v_i$ and placing it on a child of $v_i$,
{\em i.e.} a vertex $v_{i+1}$ such that $(v_i, v_{i+1}) \in E$.  ${\mathcal B}$ makes a move from $v_1$ to $v_2$ 
and gives the hand to ${\mathcal A}$, which has now to play from $v_2$, and so on. 
The first player unable to move loses the game. 
One of the two players has a winning strategy (as this game is finite in a digraph $D$ without circuit, for circuits
one extends the rules by saying that the game is lost for the player who replays a position previously reached).
Von Neumann and Morgenstern~\cite{NeMo80} proved that any directed acyclic graph  has a unique kernel, 
which is the set of winning positions for $\mathcal A$ (${\mathcal A}$ always forces
${\mathcal B}$ to play outside the kernel, until ${\mathcal B}$  cannot play anymore). 
Richardson~\cite{ri} proved later that every digraph without odd circuit has a kernel~\cite{chaty2,chaty1}. 
Berge wrote a chapter on kernels in~\cite{ber}. 
Furthermore, there is a strong connection between perfect graphs and kernels 
(see the Berge and Duchet survey~\cite{bd}). 
Some natural variants of this property are studied in various logic for 
Intelligence Artificial, some of them are definable in default
logic~\cite{cl} and used for argumentation in multi-agents  systems, 
kernels appear there as sets of coherent arguments~\cite{CaDoMe03,Do02}. 

Fernandez de la Vega~\cite{lalo} and Tomescu~\cite{to} proved 
independently that dense random digraphs with $n$ vertices and $m=\Theta ( n^2 )$ edges, 
have asymptotically almost surely a kernel. In addition, they get the few possible sizes
of a kernel and a precise estimation of the numbers of kernels.

Few years ago a new interest for these studies arises by their applications
in finite model theory. Indeed variants of kernel are 
the best properties to provide counterexamples of $0$-$1$ laws 
in fragments of monadic second-order logic~\cite{lba,lbb}. 
Goranko and Kapron showed in~\cite{gk} that such a variant 
is expressible in modal logic over almost all finite frames for frame satisfiability; 
recently Le Bars proved in~\cite{lbc} that the $0$-$1$ law fails for this logic. 

The existence of a kernel in a digraph 
has been shown NP-complete, even if one restricts this question to planar graphs 
with in- and out-degree $\leq 2$ and degree $\leq 3$~\cite{chv,cr,fr81}. 
It is somehow related to finding a maximum clique in
graphs~\cite{bbpp,lba},
which is known to be difficult for random dense graphs.

In this article, we use some generating function techniques 
to give some new results on Nim-type games played on directed graphs (or, equivalently, some new informations
on kernel of digraphs). More precisely, we deal with a family of planar digraphs 
with at most one circuit or one cycle and
we give enumerative (Theorems~\ref{theo1}, \ref{theo2}, \ref{theo3}, \ref{theo4} in Section~\ref{Sec4}) 
and asymptotics results (Theorems~\ref{theo5}, \ref{theo6}, \ref{theo7}, \ref{theo8} in Section~\ref{Sec5}) on the size of the kernel,
 the probability of winning on trees for the first player...

\section{Definitions}\label{Sec2}

We give below 
more precise definitions, readers familiar with digraphs can skip them.

Let $D = (V, E)$ be a digraph. For each $v \in V$, let 
$v^+ = \{ w \in V / (v, w) \in E\}$
and $v^- = \{ w \in V/ (w, v) \in E\}$, 
$|v^+|$ is the {\it out degree} of $v$ and 
$|v^-|$ is the {\it in degree} of $v$. 

A vertex with an in degree of 0 is called a {\em source}
(since one can only leave it) and a vertex with an out degree of 0 
is called a {\em sink} (since one cannot leave it). 
Let $U \subset V$, $U^+ = \cup _{v \in U} v^+$ and 
$U^- = \cup _{v \in U} v^-$, we denote by $D(U)$ the subgraph induced 
by the vertices of $U$.

There is a {\em path} from vertex $v$ to $w$ means that
there exists a sequence $(v_1,\ldots, v_k)$ such that 
$v_1 = v$, $v_k = w$ and $v_i \in v_{i+1}^+ \cup v_{i+1}^- $, for $i = 1\ldots k-1$. 
There is a {\em directed path} from vertex $v$ to $w$ means that
there exists a sequence $(v_1,\ldots, v_k)$ such that 
$v_1 = v$, $v_k = w$ and $v_i \in v_{i+1}^+$, for $i = 1\ldots k-1$. 

A {\it cycle } is a path $(v_1,\ldots, v_k)$ such that  $v_1 = v_k$. 
A {\it circuit} is a directed path $(v_1,\ldots, v_k)$ such that  $v_1
= v_k$. 

If $D$ contains a directed path from vertex $v$ to $w$ 
then $v$ is an {\em ancestor} of $w$ and $w$ is a {\em descendant} of
$v$. If this directed path is of length 1, then the ancestor $v$ of $w$ is
also called a {\em parent} of $w$, and $v$ a {\em child} of $w$.

$D$ is strongly connected if for each pair of vertices,  
each one is an ancestor of the other. 
$D(U)$ is a strongly connected component of $D$ if it is a maximal strongly connected subgraph  of $D$. 

$U$ is an {\it independent set} when $U \cap U^+ = \emptyset$ 
and a {\it dominating set} when $v^+ \cap U \neq \emptyset$  for any $v \in V \setminus U$. 
$U$ is a kernel if it is an independent dominating set. 

$D$ is a DAG if it is a directed digraph without circuit 
(the terminology ``directed {\em acyclic} graph'' being popular, 
we use the acronym {\em DAG}  although 
it should stands for ``directed {\em acircuit} graph'', according to the above definitions of cycles and circuits).

\newpage

\section{How to find the kernel of a digraph}\label{Sec3}

Consider digraphs satisfying the following rules:
\begin{itemize}
\item each vertex is colored either in red or in green,
\item each green vertex has at least a red child,
\item no red vertex has a red child.
\end{itemize}
It is immediate to see  that a digraph satisfying such coloring constraints
possesses a kernel, which is exactly the set of its red vertices.
It is now easy to see, e.g., 
that the circuit of length 3 has no kernel,
that the circuit of length 4 has 2 kernels, 
that any DAG has exactly one kernel.
For this last point, assume that $D$ is a DAG (directed acircuit graph). 
Algorithm~\ref{dag} (below) returns its unique kernel. It begins to color the
sinks in red and then goes up toward sources, as it is deterministic 
and as it colors at least a new vertex at each iteration, this proves
that each DAG has a single kernel. Such an algorithm was already
considered by Zermelo while studying chessgame.
\begin{algorithm}
\scriptsize{
\caption{The kernel of a DAG}
\label{dag}
\begin{algorithmic}
\STATE Input: a DAG $D=(V,E)$, Noncolored= $V$ ({\em i.e.} no vertex is colored for yet)
\STATE Output: the DAG, with all its vertices colored, the red vertices being its kernel
\WHILE{it remains some non colored vertices (Noncolored $\neq \emptyset$)}
\FORALL{$v \in$ Noncolored} 
\IF{$v$ is a sink or if all the children of $v$ are green}
\STATE color $v$ in red
\STATE color all the parents of $v$ in green
\STATE remove the colored vertices from Noncolored
\ENDIF
\ENDFOR
\ENDWHILE 
\end{algorithmic}
}
\end{algorithm}

For sure, it is possible to improve this algorithm by using the poset structure of a DAG,
and thus replacing the ``for all $v \in$ Noncolored'' line 
 by something like ``for all $v \in$ Tocolornow'' where Tocolornow is a set of candidates much smaller 
than Noncolored.

More generally, in order to color a digraph which is not a DAG, simply
split it in $p$ components which are DAGs.
Then, apply the above algorithm on each of these DAGs
(excepted the cut points that you arbitrarily fix to be red or green).
It finally remains to check the global coherence of these colorings. 
As one has $p$ cutting points (which can also be seen as $p$ branching points 
in a backtracking version of this algorithm), this leads to at most $2^p$ kernels.
This also suggests why this problem is NP: for large (dense) graph, one should need to cut at least $p\sim n$ points,
which leads to a $2^n$ complexity (lower bound). 

\begin{figure}[htbp]
\begin{center}
\includegraphics[scale=.62]{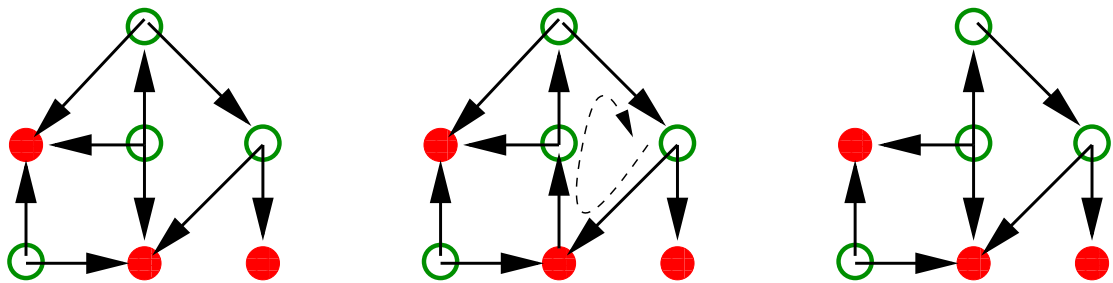}
\caption{The first digraph is a well-colored DAG (it has several cycles, but no circuit).
The second digraph is a well-colored digraph (it is not a DAG, as it contains  one circuit).
The third digraph is a DAG, but is not well colored  (the top green vertex misses a red child).
[For people who are reading a black \& white version of this article, 
red vertices are fulfilled and green vertices are empty circles.]}
\label{fig:2}
\end{center}
\end{figure}

\newpage
\section{Generating functions of well-colored unicircuit digraphs}\label{Sec4}

\def\Tru{{ T_r^\uparrow}} 
\def\Trd{{ T_r^\downarrow}} 
\def\Tgu{{ T_g^\uparrow}} 
\def\Tgd{{ T_g^\downarrow}} 
\def\Tgr{{ T_{g_r}^\uparrow}} 
\def\Set {{\rm Set}}
\def\Cyc {{\rm Cyc}}

There exists in the literature some noteworthy results on {\em digraphs} using generating functions 
(related e.g. to EGF of acyclic digraphs~\cite{Gessel,Robinson},  Cayley graphs~\cite{mishna}, (0,1)~matrices~\cite{McKay}, Erd\H os--R\'enyi random digraph model~\cite{lb}),
but as fas as we know  we give here the first example of application
to the {\em kernel} problem.

The coloring constraints mentioned in Section~\ref{Sec3} are ``local'': they are defined only in function
of each vertex and its neighbors. One nice consequence of this ``local'' viewpoint of kernels
is that it opens up a whole range of possibilities for a kind of context-free grammar approach. 
Indeed if one considers rooted labeled directed trees that are well-colored
({\em i.e.} which possesses a kernel), one can describe/enumerate them 
with the help of the five following families of combinatorial structures
(all of them being rooted labeled directed trees):
\begin{itemize}
\item $T$: all the trees with the coloring constraint
\item $\Tru$: well-colored trees with a red root (with an additional out-edge)
\item $\Trd$: well-colored trees with a red root (with an additional in-edge)
\item $\Tgu$: well-colored trees with a green root (with an additional out-edge)
\item $\Tgd$: well-colored trees with a green root (with an additional in-edge)
\item $\Tgr$: well-colored trees with a green root (with an additional
out-edge which has to be attached to a red vertex)
\end{itemize}

Those families are related by the following rules:
$$\begin{cases}
T=\Tgu \cup \Tru \\
\Tgu=g^\uparrow   \times \Set_{\geq 1}(\Tru) \times \Set(\Trd \cup  \Tgd \cup \Tgu) \\
\Tgd=g^\downarrow \times \Set_{\geq 1}(\Tru) \times \Set(\Trd \cup  \Tgd \cup \Tgu)  \\
\Tru=r^\uparrow \times \Set(\Tgd \cup \Tgr) \\
\Trd=r^\downarrow \times \Set(\Tgd \cup \Tgr)  \\
\Tgr=g^\uparrow \times \Set(\Tru \cup \Trd \cup \Tgd \cup \Tgu) 
\end{cases}$$
The {\Set } operator reflects the fact that one considers non planar trees,
{\em i.e.} the relative order of the subtrees attached to a given vertex does not matter.
The notation $\Set_{\geq 1}$ means one considers non empty set only.

\begin{figure}[htbp]
\begin{center}
\includegraphics[scale=.3]{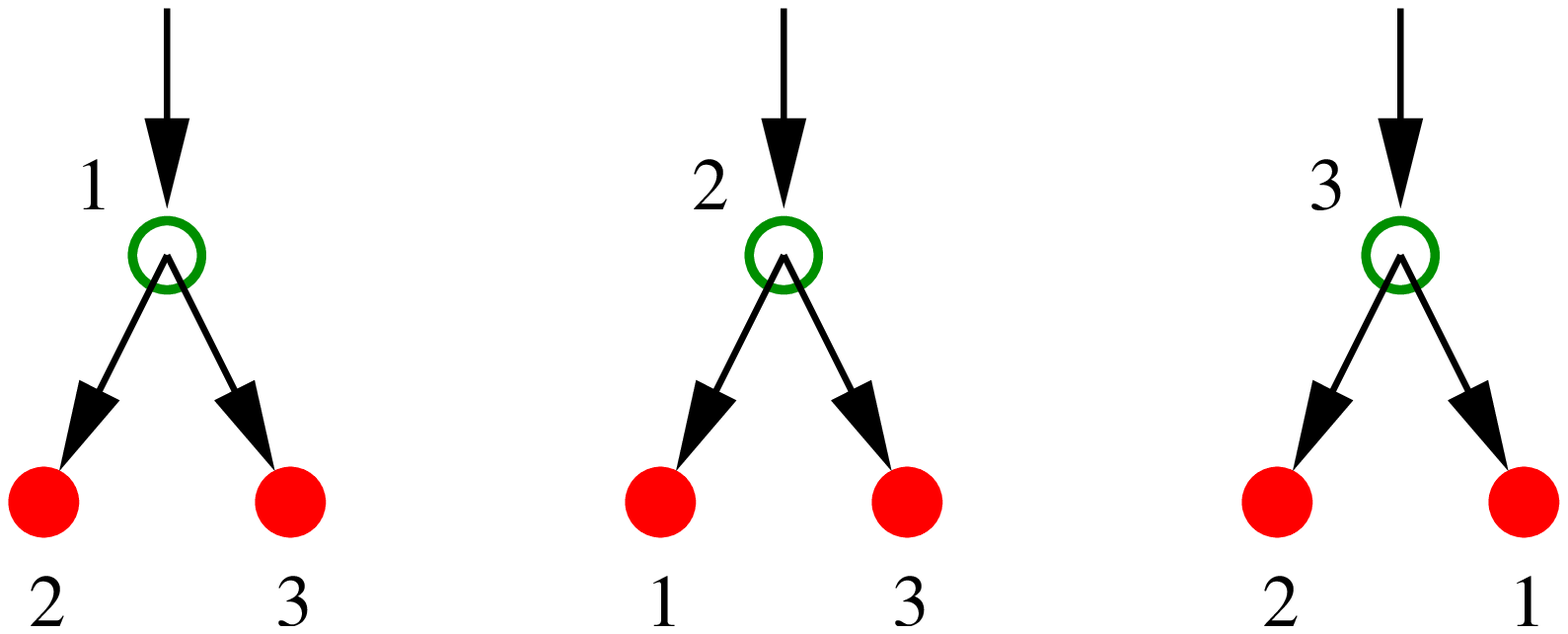}
\caption{A tree $\in \Tgd$ of size 3  and all its possible
labellings. $\Tgd$ stands for directed trees with a green root with
an additional in-edge on this root.
 \label{fig:1}}
\end{center}
\end{figure}

As we are dealing with labeled objects (we refer to Figure~\ref{fig:1} 
for the different labellings of a rooted directed tree), it is more convenient 
to use exponential generating functions, the above rules are then translated (see e.g.~\cite{hp73, FS+} 
for a general presentation of this theory of ``graphical enumeration/symbolic combinatorics'' )
into the following set of functional equations (where $z$ marks the vertices):

$$\begin{cases}
T(z)=\Tgu(z) + \Tru(z)\, ,\\
\Tgu(z)=\Tgd(z)=z (\exp(\Tru(z))-1)  \exp(\Trd(z) + \Tgd (z) + \Tgu(z))\, , \\
\Tru(z)=\Trd(z)=z \exp(\Tgd(z)+ \Tgr(z))\,. 
  \end{cases}$$

Note that $\Tgr=T$ as one has the trivial bijection  ``$\Tgr$ trees with a root without red child'' = ``$\Tru$ trees''
and ``$\Tgr$ trees with a root with at least a red child'' = ``$\Tgu$
trees''. Define now 
$T_g(z):=\Tgu(z)$ and $T_r(z):=\Tru(z)$, the above system simplifies
to:

$$\begin{cases}
T(z)=T_g(z) + T_r(z)=\Tgr(z)\, ,\\
T_g(z)=z \exp(2T(z)) -z \exp(T(z)+T_g(z))\, , \\
T_r(z)=z \exp(T_g(z)+ T(z))=T(z) \exp(-T_r(z))\,. 
  \end{cases}$$

This system has a unique solution, as the relations can be considered
as fixed point equations for power series.
Their Taylor expansions are:
$$T(z)=   z+4\frac{z^2}{2!}+36\frac{z^3}{3!}+512\frac{z^4}{4!}+10000\frac{z^5}{5!}+248832\frac{z^6}{6!}+7529536\frac{z^7}{7!}+O(z^8)\,,$$
$$T_g(z)=\quad 2\frac{z^2}{2!}+15\frac{z^3}{3!}+232\frac{z^4}{4!}+4535\frac{z^5}{5!}+114276\frac{z^6}{6!}+3478083\frac{z^7}{7!}+O(z^8)\,,$$
$$T_r(z)=z+2\frac{z^2}{2!}+21\frac{z^3}{3!}+280\frac{z^4}{4!}+5465\frac{z^5}{5!}+134556\frac{z^6}{6!}+4051453\frac{z^7}{7!}+O(z^8)\,.$$

For sure, the $i$-th coefficients of these series are divisible by $i$, as we are dealing with rooted object.
Here are the 3 generating functions of the corresponding unrooted trees:
$$T^{unr.}(z)=   z+2\frac{z^2}{2!}+12\frac{z^3}{3!}+128\frac{z^4}{4!}+2000\frac{z^5}{5!}+41472\frac{z^6}{6!}+1075648\frac{z^7}{7!}+O(z^8)\,,$$
$$T_g^{unr.}(z)=\quad \frac{z^2}{2!}+5\frac{z^3}{3!}+58\frac{z^4}{4!}+907\frac{z^5}{5!}+19046\frac{z^6}{6!}+496869\frac{z^7}{7!}+O(z^8)\,,$$
$$T_r^{unr.}(z)=z+\frac{z^2}{2!}+7\frac{z^3}{3!}+70\frac{z^4}{4!}+1093\frac{z^5}{5!}+22426\frac{z^6}{6!}+578779\frac{z^7}{7!}+O(z^8)\,.$$

\begin{figure}[htbp]
\begin{center}
\includegraphics[scale=.52,angle=-90]{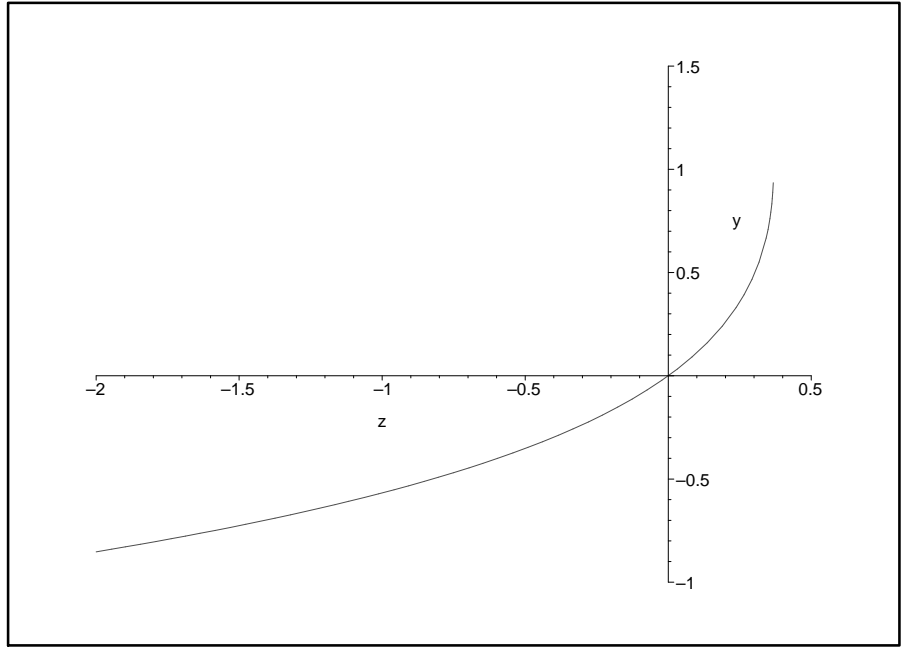}
\caption{The Cayley tree function $C(z)$ 
goes from $-\infty$ for $z\sim-\infty$ to~1 at $z=\frac{1}{e}$. It satisfies $C(z)=z \exp(C(z))$.}
\label{fig:4}
\end{center}
\end{figure}

Of course, trees are DAG and therefore have a unique kernel. This implies that $T(z)$ is exactly 
the exponential generating function of directed rooted trees, {\em i.e.}
$$T(z)=C(2z)/2  \text{ and } T_n=(2n)^{n-1}$$
where $C(z)$ is the Cayley function 
(see Figure~\ref{fig:4} and references~\cite{Ca89,CoJeKn97}), defined by 
$$C(z)=z \exp(C(z))=\sum_{n\geq 1} n^{n-1} \frac{z^n}{n!}\,.$$

Solving the set of equations for $T, T_g$ and $T_r$  finally leads to 
\begin{theorem}[Enumeration of well-colored trees]\label{theo1}\quad\newline 
By ditrees, we mean well-colored rooted labeled directed trees. By well-colored, we mean each green vertex has at least a red child, each red vertex has no red child.\\
The exponential generating function of ditrees is given by $T(z)=C(2z)/2$, \\ 
the EGF of ditrees with a red root is given by $T_r(z)=-C(-C(2z)/2)$, \\ 
the EGF of ditrees with a green root is given by
$T_g(z)=C(2z)/2+C(-C(2z)/2)$, \\ 
where $C(z)$ is the Cayley tree function  $C(z)=z \exp(C(z))$.\\
The EGF for the unrooted equivalent objects can be expressed in terms of the rooted ones:
$$T^{unr.}=T-T^2\,, \qquad T_g^{unr.}=T^{unr.}-T_r^{unr.}\,, \text{\quad and  \quad} T_r^{unr.}=2T-2T T_r+T_r-2T/T_r+T_r^2/2\,.$$
\end{theorem}
\begin{proof}
The formulae for $T, T_r$ and $T_g$ can be checked using the definition
of $C(z)$ in the fix-point equations in the simplified system above.
The fact that the GF for unrooted trees  can be expressed in terms of the GF of rooted ones
can be proven by integration of the Cayley function, or by a
combinatorial splitting argument on trees.
\end{proof}

We can go on and enumerate the different possibilities of circuits
for a well-colored digraph. 
They can be described as
\[ \Cyc( g ) \quad \cup \quad  \Cyc( r {\text \small \rightarrow}  \{g
{\text \small \rightarrow}\}^+ ) \]
This reflects the fact that either a circuit is made up of green vertices only, 
or it contains some red vertices, but they have to be followed by at least a green vertex.
NB: Whether one counts or not the cycles of length 1 ({\em i.e.} a single red or green vertex) 
will only modify the first term of the generating function.
Symbolic combinatorics~\cite{FS+} translates the above cycle decompositions in the following 
function:
 $$\ln\left(\frac{1}{1-g}\right) + \ln\left(\frac{1}{1-\frac{rg}{1-g}}\right) $$
where $r$/$g$ mark the number of red/green vertices.
This leads to the following Theorem:

\begin{theorem}[Enumeration of possible well-colored circuits]\quad\newline\label{theo2}
The exponential generating function of possible well-colored  circuits is given by  
$$L(z)=-\ln(1-z-z^2)=
z+3\frac{z^2}{2!}+8\frac{z^3}{3!}+42\frac{z^4}{4!}+264\frac{z^5}{5!}+2160\frac{z^6}{6!}+20880\frac{z^7}{7!}+O(z^8)\,. $$ 
Its coefficients satisfy $L_n = (n-1)!\, \left(\phi^n+\left(1-\phi \right)^n\right)$, 
where $L_n$ are known as the $n$-th Lucas number (usually defined by the recurrence $L_{n+2}=L_{n+1}+L_n, L_1=1, L_2=3$) 
and where $\phi=(1+ \sqrt 5)/2$ is the golden ratio.
\end{theorem}
Note that a reverse engineering lecture of this generating function leads to the simpler
decomposition $\Cyc(g \cup rg)$, which also explains the recurrence!
Now, the following decomposition of possible {\em cycles} is trivially
related to the decomposition of possible {\em circuits}:
$$\Cyc( r \times \{{\text \small \rightarrow}  g  \cup {\text \small
\leftarrow} g\}^+ \times  \{{\text \small \rightarrow} 
\cup {\text \small \leftarrow} \}  )   \quad \cup \quad  \Cyc( g
{\text \small \rightarrow}  \cup g {\text \small \leftarrow} )$$
leads to the EGF $-\ln(1-2z-4z^2)$ whose coefficients are, with no surprise, $2^n L_n$.

\begin{figure}[hbtp]
\begin{center}
\includegraphics[scale=.29]{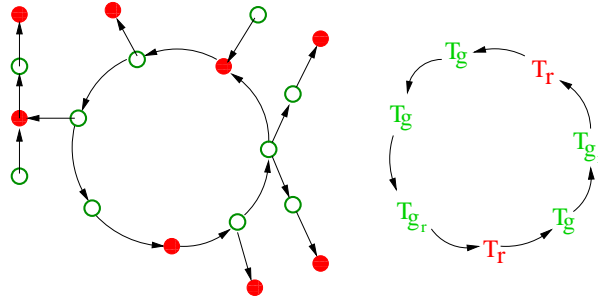}
\caption{Unicircuit digraphs consist in a circuit with attached trees on it. The left picture above
is a unicircuit digraph, to the right, we give its ``canonical decomposition'' as a circuit of atoms
which are trees. Any well-colored unicircuit digraph 
has such a ``canonical decomposition''.}\label{fig:3}
\end{center}
\end{figure}

Using the decomposition given in Figure~\ref{fig:3}, 
one obtains the generating function for unicircuits:
\begin{theorem}[Enumeration of unicircuit well-colored digraphs]\quad\newline\label{theo3}
The EGF of unicircuit well-colored digraphs is 
\begin{eqnarray*}
U(z)&=&T^{unr.}-T_g + \ln\left(\frac{1}{1- (T_g+T_{g_r}T_r) }\right) \\
&=& -\frac{C(2z)^2}{4}-C(-\frac{C(2z)}{2}) -\ln\left(1-\frac{C(2z)}{2}-C(-\frac{C(2z)}{2}) + {C\big(-\frac{C(2z)}{2}\big) \frac{C(2z)}{2}}\right)  \\
&=& z+4\frac{z^2}{2!}+30\frac{z^3}{3!}+452\frac{z^4}{4!}+8840\frac{z^5}{5!}+224832 \frac{z^6}{6!}+6909784 \frac{z^7}{7!}+O(z^8)\,,
\end{eqnarray*}
where $C(z)$ is the Cayley tree function  $C(z)=z\exp(C(z))$.
\end{theorem}

Now, consider the larger class of unicycle digraphs
(digraphs which have 0 or 1 cycle). Recall that a circuit is a cycle,
but a cycle is not necessarily  a circuit. In order to get a ``canonical
decomposition'' for unicyle digraphs (similar to the one given in
Fig.~4 for unicircuit digraphs), one considers 3 cases:

\begin{itemize} 
\item Either the graph has no cycle, those graphs are counted by $T^{unr.}$. 
\item  Either it is a cycle with  only $T_g$  trees branched on it
({\em i.e.\/} no red nodes in the cycle), those graphs
are counted by $(\ln\left(\frac{1}{1-2T_g}\right)-2Tg-4Tg^2/2)/2+Tg^2/2$, where 
$2Tg$ corresponds to  $T_g \times  \{\rightarrow \cup \leftarrow\}$,
one removes cycles of length 1 and 2 from the logarithm  (this explains the
$-2Tg-4Tg^2/2$ term) and one divides the whole formula 
by 2  because one has to take into account the fact the cycle can be read 
clockwise or not, and one adds the only legal cycle of length 2.
\item Either the graph contains a cycle with some red nodes and then one considers the following
 possible ``bricks'':
$$\begin{cases}
T_r \leftarrow  T_{g_r} \leftarrow    \\  
T_r \leftarrow  T_{g_r}  \rightarrow \qquad \text{(but not a cycle of length 2, because multiarcs are not allowed)  }    \\
T_r \rightarrow \left( T_g  \{\rightarrow \cup \leftarrow\}\right)^* T_g \leftarrow \qquad \text{(but not a cycle of length 2)} \\
T_r \rightarrow  \left( T_g  \{\rightarrow \cup \leftarrow\}\right)^*   T_{g_r} \rightarrow  \\
T_r \leftarrow  T_{g_r} \{\rightarrow \cup \leftarrow\}  \left( T_g  \{\rightarrow \cup \leftarrow\}\right)^* T_{g}  \leftarrow  \\
T_r \leftarrow  T_{g_r} \{\rightarrow \cup \leftarrow\} \left( T_g \{\rightarrow \cup \leftarrow\}\right)^* T_{g_r}  \rightarrow 
\end{cases}$$
\end{itemize}

\begin{theorem}[Enumeration of unicycle well-colored digraphs]\quad\newline\label{theo4}
The EGF of unicycle well-colored digraphs is 
\begin{eqnarray*}
V(z)&=&T^{unr.}+\frac{1}{2}\ln\left(\frac{1}{1-2T_g}\right) -Tg-Tg^2/2
-TrTg/2 -T_r T/2 \\
& &+ 
\frac{1}{2}\ln\left(\frac{1}{1- \left( 2 T_r T_{g_r}+
 \frac{T_r T_{g}+T_r T_{g_r}+ 2  T_r T_{g_r} T_g+  2  T_r T_{g_r}^2  }{1-2T_g}
\right)}\right) \\
&=&T^{unr.}-T+T_r -T^2/2 - \ln(1+T_r) - \frac{1}{2} \ln(1-2T)  \\
&=& z+4\frac{z^2}{2!}+36\frac{z^3}{3!}+692\frac{z^4}{4!}+15920\frac{z^5}{5!}+458622 \frac{z^6}{6!}+15559264 \frac{z^7}{7!}+O(z^8)\,. 
\end{eqnarray*}
where $T$, $T_g$, $T_r$, and $T^{unr.}$  are given in Theorem~\ref{theo1}.
\end{theorem}

Note that in the two theorems above, any given non-colored graph 
is counted with multiplicity 0, 1 or  2 (if there are 0, 1 or 2 ways
to color it). We explained in Section 3 that a multiplicity larger
than 2 was not possible for unicycle digraphs. We enumerate in the
following proposition those  with exactly 2 possible colorations.
\begin{prop}[Enumeration of unicycle digraphs with two kernels]\quad\newline\label{prop1}
The EGF of unicycle digraphs with 2 kernels is 
$$
D(z)=-\ln\sqrt{1+C(-C(2z)/2)^2}\,,$$
where $C(z)$ is the Cayley tree function $C(z)=z\exp(C(z))$.
\end{prop}
{\em Remark:} From the definition of cycle/circuit, $D(z)$ is also the
EGF of {\em unicircuit} digraphs with 2 kernels.
\begin{proof}
Let $\mathcal D$ be the set of unicycle digraphs with 2 kernels.
First, it is easy to see that $\Cyc(T_r^2) \subset \mathcal D$ (with a
sligth abuse of notation, as we first  only consider the shape, not
the coloration of the $T_r$ trees): 
simply color the nodes in the cycle alternatively in green and red, 
and switch the colors of a part of attached trees, if needs be.
 
We now prove the next step $ \mathcal D \subset \Cyc(T_r^2)$:
Take a unicycle graph in $\mathcal D$, it means at least one of its 
vertex can be colored both green and red. Such a vertex $v$ 
can be taken, without loss of generality, in the circuit (from the
above remark, the cycle is in fact a circuit).
 [If it were not the case, all bi-colorabled vertices would be in the tree components, but then one
could split our graph to get DAGs which are known to be uniquely colorable].  
But when $v$ is red, it implies it has only $T_g$ trees attached to it,
which means than when it gets green, the next node in the circuit has
be red (and was previously green!). This implies alternation red/green
(and even length for parity reasons) for all the nodes in the circuit.

This leads to a canonical decomposition 
$$\Cyc( T_r^2 )\,.$$
If one divides by 2 for the (anti)clockwise readings, this leads to
the Theorem.
\end{proof}

\bigskip
Most of these results (and also the computations of Section~\ref{Sec5}
hereafter) were checked with the computer algebra system Maple.
A worksheet corresponding to this article 
is available at  {\tt
http://algo.inria.fr/banderier/Paper/kernels.mws}
(or  {\tt kernels.html}), it uses 
the Algolib librairy,
downloadable at {\tt http://algo.inria.fr/libraries/}).

\section{Asymptotics}\label{Sec5}

In this section, we give asymptotic results for $n \rightarrow +\infty$.

\begin{theorem}[Proportion of trees with a green/red root]\quad\newline \label{theo5}
Asymptotically $\frac{1-\lambda}{1+\lambda}\approx 47.95\%$ of the
trees have a green root,
where the constant $\lambda\approx 0.351733$ is defined as the unique real root of $2\lambda= \exp(-\lambda)$.

A more pleasant way to formulate this Theorem consists in considering 
Nim-type games (first player who cannot move loses) 
on directed trees where the tree and the starting position are chosen uniformly at random.
The strategies of the two players being optimal, 
the first player has then a probability of 47.95\% (asymptotically) to win the game.
(Recall that if the starting position can be chosen by the first player, then he will always win.)
\end{theorem}
\begin{proof}
The key step of this result and the following ones are the following
expansions (derived from the expansion of the Cayley function) for $T$, $T_r$ and $T_g$:

$$T(z)  \sim \frac{5}{6} - \frac{1}{\sqrt 2} \sqrt{1-2ez} +O(1-2ez)$$
$$T_r(z)\sim \lambda - \frac{\lambda \sqrt 2}{1+\lambda} \sqrt{1-2ez} +O(1-2ez)$$
$$T_g(z)\sim \frac{1}{2}-\lambda - \frac{1}{\sqrt 2} \frac{1-\lambda }{1+\lambda} \sqrt{1-2ez} +O(1-2ez)\,,$$
where the constant $\lambda$ is defined as
$\lambda:=T_r(\frac{1}{2e})\approx 0.351733$.

By Pringsheim theorem~\cite{FS+},  as $T_r(z)$ has nonnegative
coefficients, then $T_r(z)$  has a positive singularity.
As coefficients of $T_r$ are smaller
than coefficients of $T$, 
its radius of convergence belongs to  $[0,1/(2e)]$. 
Now, $-C(2z)/2$ is negative on this interval, 
and thus $C(-C(2z)/2)$ is analytic there, and $1/(2e)$ is therefore
its only possible dominating singularity.
The radius of $T_g$ follows from $T=T_r+T_g$. 
The theorem follows by considering  $\frac{[z^n] T_g(z)}{[z^n] T(z)}=\frac{1-\lambda}{1+\lambda}-\frac{\lambda^2(\lambda+4)}{(1+\lambda)^5}\frac{1}{n}+O(\frac{1}{n^2})$.
\end{proof}

\begin{theorem}[Proportion of red vertices in possible circuits]\quad\newline \label{theo6}
Asymptotically $\frac{1}{2}-\frac{1}{2\sqrt 5}\approx 27.63\%$ of the vertices of a possible circuits are red.
\end{theorem}
\begin{proof}
One has to considerer the following bivariate generating function
(exponential in $z$, ordinary in $u$): $\ln\left(\frac{1}{1-(z+uz^2)}\right)$.
The wanted proportion is then given by $\frac{[z^n] \partial_u
F(z,1)}{[z^n] F(z,1)}$,
where [$z^n] \partial_u F(z,1)$ means the $n$-th coefficient of ``the
derivative with respect to $u$ of $F(z,u)$, then evaluated at $u=1$''. 
\end{proof}

Then, one can wonder if the asymptotic density of  well-colored
unicircuit graphs is more than 50\% or even if it is 100\%?  The
following theorem gives the answer:

\begin{theorem}[Proportion of well-colored unicircuit digraphs]\quad\newline \label{theo7}
The  proportion of well-colored graphs amongst unicircuit digraphs
is asymptotically:
$$\frac{3 \lambda^3+\lambda^2-\lambda-1}{(1+\lambda)^2(\lambda-1))} \approx 92.65\%$$
 where $\lambda$ is the constant defined in Theorem~\ref{theo5}.
\end{theorem}
\begin{proof}
Relies on a singularity analysis of the generating function of Theorem~\ref{theo3},
with the expansions given in Theorem~\ref{theo5}.
Note that some unicircuit digraphs can have 2 kernels,
so one has to perform the following asymptotic expansions:
$$\frac{[z^n] U(z)-D(z)}{[z^n]F(z)}\approx 92.65 -\frac{0.12}{n}+O(\frac{1}{n^2})\,, $$
where $F(z)=T^{unr}(z)+\ln(\frac{1}{1-T(z)})-T(z)$
is the EGF of (non-colored) unicircuit digraphs.
\end{proof}

For sure, it one considers now the asymptotic density of well-colored
{\em  unicircuit} graphs, the proportion should be larger, as one only
adds DAGs (which are all well-colorable). The following theorem gives the
noteworthy result that unicircuit graphs are in fact almost surely
well-colored:

\begin{theorem}[Proportion of well-colored unicycle digraphs]\quad\newline \label{theo8}
There is asymptotically a proportion of 
$1- \frac{2\lambda^3 \sqrt 2}{(1+\lambda)^2 (1-\lambda) \sqrt {\pi }} \frac{1}{\sqrt n} \approx
1-\frac{0.05}{\sqrt n}$  of well-colored graphs amongst unicycle digraphs of size $n$,
where $\lambda$ is the constant defined in Theorem~\ref{theo5}.
\end{theorem}
\begin{proof}Relies on a singularity analysis of the generating function of Theorem~\ref{theo4},
with the expansions given in Theorem~\ref{theo5}.
Note that some unicycle digraphs can have 2 kernels,
so one has to consider $$\frac{[z^n] V(z)-D(z)}{[z^n]G(z)}\,,$$
where $G(z)=T^{unr}(z)+\frac{1}{2}\ln(\frac{1}{1-2T(z)})-T(z)-T(z)^2/2$
is the EGF of (non-colored) unicycle digraphs (one substracts $T^2/2$
because amongst the 4 graphs with a cycle of length 2 created 
by the $\ln(\frac{1}{1-2T(z)})$ part, 3 are not legal: 
1 was already counted because of
symmetries, and the other 2 have in fact a multiple arc, 
whereas it is forbidden for our digraphs).
\end{proof}

Finally, if one considers graphs with at most $k$ cycles, 
it means one has more cutting points, which  relaxes the constraints for
well-colarability (=existence of kernel). According to the above
results, this implies an asymptotic density of one. 
This gives as a corolary of our results, that all these families
have almost surely a kernel. A kind of ``limit case'' is dense graphs,
for which  some results already
mentionned by  Fernandez de la Vega~\cite{lalo} and Tomescu~\cite{to}
give that they have indeed almost surely a kernel.

\newpage
\section{Conclusion}\label{Sec6}

It is quite pleasant that our generating function approach allows
to get new results on the kernel problem, giving {\em e.g.} the proportion of graphs satisfying a given property,
and new informations on Nim-type games for some families of graphs.

As a first extension of our work, it is possible to apply classical
techniques from analytic combinatorics~\cite{FS+} in order to get 
informations on standard deviation, higher moments, and 
limit laws  for statistics studied in Section~\ref{Sec5}.

Another extension is to get closed form formulas 
 for bicircuit/bicycles digraphs, (the generating function  involves
the derivative of the product of two logs and the asymptotics are performed like in our
Section~\ref{Sec5}). It is still possible
(for sure with the help of a computer algebra system)
to do it for 3 or 4 cycles but the ``canonical 
decompositions'' and the  computations get cumbersome.
In order to go on our analysis far beyond low-cyclic graphs,
one needs an equation similar to the one given by E.M. Wright~\cite{Wr77,Wr80} for graphs.
Let $\FW_{\ell}$ be the family of well-colored digraphs with $\ell$ edges more than vertices, ($\ell \geq -1$). 
It is possible to get an equation ``\`a la Wright'' for $\FW_{\ell}$ 
by pointing any edge (except edges linking a green vertex to a red one) in a well-colored digraph.
It is however not clear for yet if and how such equations can be simplified in order to get a recurrence
as ``simple/nice'' to the one that Wright got for graphs, thus opening the door to asymptotics
and threshold analysis  beyond the unicyclic case.

\vskip 5mm
{\bf Acknowledgements.}
This work was partially realized while the second author (Jean-Marie Le Bars) benefited 
from a delegation CNRS in the Laboratoire d'Informatique de Paris Nord (LIPN).
We would like to thank Guy Chaty  and Fran{\c c}ois L\'evy 
for a friendly discussion and their bibliographical references~\cite{chaty2,chaty1}.
We also thank an anomymous referee for detailed comments on a somehow
preliminary version of this work.

\bibliographystyle{plain}
\bibliography{fpsac04}
\vskip 1cm
\end{document}